\documentclass{article}
\usepackage[utf8]{inputenc}
\usepackage{amssymb,latexsym,amsmath} 
\usepackage{indentfirst}
\usepackage{algorithm}
\usepackage{algpseudocode}
\usepackage{graphicx}
\graphicspath{ {Desktop/} }

\usepackage{graphicx}

\begin{document}

\title{Improvement of the Cascadic Multigrid Algorithm with a Gauss Seidel Smoother to Efficiently Compute the Fiedler Vector of a Graph Laplacian}
\author{Shivam Gandhi - Tufts University Department of Mathematics \thanks{shivam.jgandhi@gmail.com}}
\date{November 2015}
\maketitle

\section{Abstract}
In this paper, we detail the improvement of the Cascadic Multigrid algorithm with the addition of the Gauss Seidel algorithm in order to compute the Fiedler vector of a graph Laplacian, which is the eigenvector corresponding to the second smallest eigenvalue. This vector has been found to have applications in graph partitioning, particularly in the spectral clustering algorithm. The algorithm is algebraic and employs heavy edge coarsening, which was developed for the first cascadic multigrid algorithm. We present numerical tests that test the algorithm against a variety of matrices of different size and properties. We then test the algorithm on a range of square matrices with uniform properties in order to prove the linear complexity of the algorithm. 

\section{Introduction}

The Fiedler vector has seen numerous applications within computational mathematics, primarily within the fields of graph partitioning and graph drawing [1]. In particular, we require eigenvalues and eigenvectors for a successful run of the spectral clustering algorithm that partitions a network into clusters[2]. Although many languages have built in eigenvalue methods, the spectral clustering method requires a specialized eigenvalue algorithms to account for massive network size. In fact, spectral clustering becomes unfeasible for networks of size over 1000 as computing the eigenvalues through matrix inversion becomes inefficient. Therefore, we require specialized multigrid algorithms to find the eigenvectors and eigenvalues in less than $O(N^3)$ time.
    
The Cascadic Multigrid Algorithm is an effective method for computing the second largest eigenvalue and eigenvector, where the eigenvector is called the Fiedler vector. The main methods for calculating eigenvalues and eigenvectors include the Lanczos method and the power method. However, these methods become unfeasible for large matrices ($|V| > 1000$). Furthermore, many networks can have over $1000$ nodes which correlates to a matrix with dimension higher than $1000$. For this reason, we require the cascadic multigrid algorithm, as it solves the eigenvalue problem on coarser levels and projects the solution upwards until the solution is projected to the original matrix [3].
    
When calculating the eigenvalues and eigenvectors of symmetric positive definite matrices more generally, the Jacobi and PCG methods provide good approximations. These can be extended to calculating the Fiedler vector and its eigenvalue [1].

In this paper, we improve upon the previously made Cascadic Multigrid Algorithm by introducing a Gauss Seidel smoother on each level. We employ the previously established heavy edge coarsening that selects the edge with the heaviest weight between two vertices. The refinement procedure continues to use power iteration on a modified matrix. This method does not require the inversion of matrices, unlike Rayleigh-Quotient iteration, thereby making it a much more optimal method. As always, the eigenvector calculated on coarse levels are projected to a finer level with interpolation matrices. The eigenvector that gets calculated and projected up to a higher level serves as the guess for the next Gauss Seidel iteration on the finer level. Finally, on the highest level, we calculate the Rayleigh Quotient to achieve the eigenvalue.
    
The paper is organized to provide a logical introduction to the algorithm. In section 3, we provide definitions and background knowledge required to understand multigrid methods. We also introduce heavy edge coarsening, power iteration, and the Gauss Seidel method. This culminates in a presentation of the algorithm developed. In section 4, we introduce the numerical tests of the algorithm. We compare the algorithm to the previous cascadic multigrid algorithm and various other multigrid algorithms that are meant to calculate the Fiedler vector. We also compare spectral clustering with the build in eigenvalue calculating function in MATLAB to spectral clustering that employs our algorithm to show efficient. In the final section, we wrap up the paper and discuss future improvements to multigrid algorithms that calculate Fiedler vectors.

\section{Modified Cascadic MG Method for Computing the Fiedler Vector}

First we formally introduce the concepts of the graph Laplacian and the Fiedler vector. A weighted graph $G = (V,E,w)$ is undirected if the edges are unoriented. 

\textbf{Definition 2.1} $G = (V,E,w)$ is a weighted graph. The Laplacian of $G$, $L(G) \in \mathbb{R}^{n \times n}$, shortened to $L$, where $n = |V|$, is denoted as follows

\begin{displaymath}
   L(G)_{(i,j)} = \left\{
     \begin{array}{lr}
       d_{v_i} &, $ if $ i = j\\
       -w_{(i,j)} &, $ if $ i \neq j\\
     \end{array}
   \right.
\end{displaymath} 
where $d_{v_i}$ is the degree of vertex i and $w_{i,j}$ is the weight of the edge connecting $v_i$ and $v_j$.

This Laplacian is positive semi-definite and diagonally dominant, and the sum of and row or column of L is zero. This makes the smallest eigenvalue 0 with the corresponding vector $[1,1,...,1]^T$. We are particularly interested in the second smallest eigenvalue and eigenvector.
    
\textbf{Definition 2.2} The second smallest eigenvalue of the Laplacian of a graph $G$ is called the algebraic connectivity. This eigenvalue must be greater than or equal to 0. The corresponding eigenvector $\phi_2$ is called the Fiedler vector of G.

The importance of the Fiedler vector is detailed in [4, 5].

It is important to note that the coarsest graph must be very small in size at around $|V| < 25$. A direct power iteration is used at this coarsest level to obtain an eigenvector. Afterwards, the eigenvector is projected upwards and then smoothed using Guass-Seidel. 

We now introduce heavy edge coarsening for our cascadic algorithm. In our algorithm, $L^i \in \mathbb{R}^{n_{i} \times n_{i}}$. Heavy edge coarsening is iterated on the graph Laplacian in order to create multiple levels for solving. This algorithm makes up the setup phase.

\begin{algorithm}
\caption{Heavy Edge Coarsening}
\begin{algorithmic}[1]
\Procedure {HEC}{L}
\State $c \leftarrow 0$
\State $p \leftarrow randperm(n_i)$
\State $q \leftarrow zeros(n_i,1)$
\For {$i=1 \rightarrow n_i$}
\If {$q(p(i)) = 0$} 
\State $m \leftarrow argmin(L(:,p(i)))$
\If {$q(m) = 0$}
\State $c \leftarrow c+1$
\State $q(m) = c$
\State $q(p(i)) = c$
\Else
\State $q(p(i)) = q(m)$
\EndIf
\EndIf
\EndFor
\State $I_{i}^{i+1} \leftarrow zeros(c,n_i)$
\For {$i=1 \rightarrow n_i$}
\State $I_{i}^{i+1}(q(i),i) = 1$
\EndFor
\EndProcedure
\end{algorithmic}
\end{algorithm}

Heavy edge coarsening is further detailed in [3], and several properties of the algorithm are proved as well. 

Next, we formally introduce the Gauss Seidel method. This method takes a guess vector and solves a linear system using that guess. In our algorithm, the we use the vector projected upwards from the coarser level as the guess. This was similar to power iteration, as we used the projected vector as the first guess for power iteration as well. The values A and b are the original values in the linear system $Ax = b$. X0 is our initial guess to the solution of this system. N denotes the number of iterations allowed while tol represents the tolerance of error. The algorithm outputs a solution to $Ax = b$ within our denoted error. 

\begin{algorithm}
\caption{Gauss Seidel}
\begin{algorithmic}[1]
\Procedure {G-S}{A, b, X0, tol, N}
\State $k \leftarrow 1$
\While {$k \leq N$}
\For {$i = 1 \rightarrow n$}
\State $x_i = 1/a_{ii}[- \sum\limits_{j = 1}^{i-1} (a_{ij}x_j) - \sum\limits_{j = i+1}^{n}(a_{ij}X0_j) + b_i]$
\If {$ |x - X0| < tol$}
\State output $[x_1, x_2, ..., x_n]$
\EndIf
\State $k = k+1$
\For {$i = 1 \rightarrow n$}
\State $X0_j = x_i$
\EndFor
\EndFor
\EndWhile
\State Output $[x_1, x_2, ..., x_n]$
\EndProcedure
\end{algorithmic}
\end{algorithm}

We discuss two theorems that confirm that the Gauss Seidel method will converge to a solution in our multigrid algorithm.

\textbf{Theorem 2.1}: The Gauss Seidel method converges if A is symmetric positive definite or if A is strictly or irreducibly diagonally dominant.

\textbf{Theorem 2.2}: Let A be a symmetric positive definite matrix. Then the Gauss-Seidel method converges for any arbitrary choice of initial approximation $x$.

A proof of these theorems can be found in [6]. All of our graph Laplacians on all levels are symmetric positive definite and diagonally dominant therefore the Gauss Seidel method will converge on all levels.

With our component algorithms defined and sufficiently detailed, we can now outline the procedure for our algorithm. We begin with a setup phase that has heavy edge coarsening set up the levels on which we do computations. After this, we solve the eigenvalue problem on the coarsest level. We then begin projecting our eigenvector upwards and using Gauss Seidel on finer and finer levels until we get to the finest level, our original matrix. At this level, we use Gauss Seidel one last time to yield the Fiedler vector and then calculate the Rayleigh quotient for the algebraic connectivity. We input the finest level graph Laplacian and the algorithm outputs the Fiedler vector and corresponding eigenvalue.

\begin{algorithm}
\caption{Gauss Seidel Cascadic Multigrid}
\begin{algorithmic}[1]
\Procedure {Step 1: Setup Phase}{L}
\State $i = 0$
\While {$n_i > 25$}
\State $I_i^{i+1} \leftarrow HEC(L^i)$
\State $L^{i+1} = I_{i}^{i+1}L^{i}(I_{i}^{i+1})^T$
\State $i = i+1$
\EndWhile
\State $J \leftarrow i$
\EndProcedure
\Procedure {Step 2: Coarsest Level Solving Phase}{$L^J$}
\State $ y^(J) \leftarrow GS(L^J, rand(n_J))$
\EndProcedure
\Procedure{Step 3: Cascadic Refinement Phase}{$y^J$, $L$}
\For {$j = J-1 \rightarrow 0$}
\State $y^{j} = (I_{i}^{i+1})^{T}y^{(j+1)}$
\State $y^{j} \leftarrow GS(L^j, y^{j})$
\EndFor
\EndProcedure
\end{algorithmic}
\end{algorithm}

Structurally, this algorithm is similar to other multigrid algorithms in that it begins with a setup phase and solves on the coarsest level upwards. It is nearly identical to the Cascadic Multigrid Algorithm with the sole difference being in the Gauss Seidel replacing power iteration.

\section{Numerical Tests}

We perform numerical tests on a variety of graphs listed on Table 4.1. The graphs were taken from the University of Florida Sparse Matrix Collection [7]. The computations were performed on an HP Envy with a 2.40 GHz Intel Core i7 Processor with 8.00 GB of RAM. We consider the performance of the Gauss Seidel Cascadic Multigrid Algorithm to matrices with over 8000 nodes. We use a tolerance $(u^k,u^{k-1}) > 1 - 10^{-6}$. 

\begin{tabular}{l||l|l|l}
Matrix Name & Matrix Size & Matrix Edges & CGMG runtime (s) \\ \hline
barth5      & 15606       & 45878        & 0.371467         \\ \hline
bcsstk32    & 44609       & 985046       & 1.242307         \\ \hline
bcsstk33    & 8738        & 291583       & 0.381135         \\ \hline
brack2      & 62631       & 366559       & 1.307903         \\ \hline
copter1     & 17222       & 96921        & 0.42307          \\ \hline
ct2010      & 67578       & 168176       & 1.265944         \\ \hline
halfb       & 224617      & 6081602      & 6.694857         \\ \hline
srb1        & 54924       & 1453614      & 1.582835         \\ \hline
wing\_nodal & 10937       & 75488        & 0.40845         
\end{tabular}

Next, we show that the algorithm is $O(N)$. We run the algorithm on uniform square arrays of various sizes and show that the runtime increases linearly according to the matrix size. The amount of nodes and edges increases linearly therefore we can expect the runtime of the algorithm to also increase linearly. Because multigrid algorithms run in linear time, it is important that the Gauss Seidel smoother does not change the runtime, otherwise it would be an inferior algorithm to use. The r value is very close to 1, indicating that the algorithm does in fact have $O(N)$ complexity.

\begin{tabular}{l|l}
Matrix Nodes                 & Time (seconds) \\ \hline
106276                       & 1.921614       \\ \hline
178929                       & 3.220836       \\ \hline
232324                       & 4.088426       \\ \hline
276676                       & 5.344172       \\ \hline
303601                       & 5.684314       \\ \hline
374544                       & 7.178143       \\ \hline
425104                       & 7.811554       \\ \hline
564001                       & 10.565033      \\ \hline
657721                       & 11.704087      \\ \hline
705600                       & 12.936846      \\ \hline
736164                       & 13.768696      \\ \hline
753424                       & 13.843865      \\ \hline
762129                       & 14.799933      \\ \hline
788544                       & 14.613115      \\ \hline
795664                       & 15.51262       \\ \hline
799236                       & 16.808463      \\ \hline
848241                       & 16.922279      \\ \hline
851929                       & 16.233831      \\ \hline
915849                       & 17.257426      \\ \hline
956484                       & 19.349795    
\end{tabular}

\includegraphics[scale = 0.5] {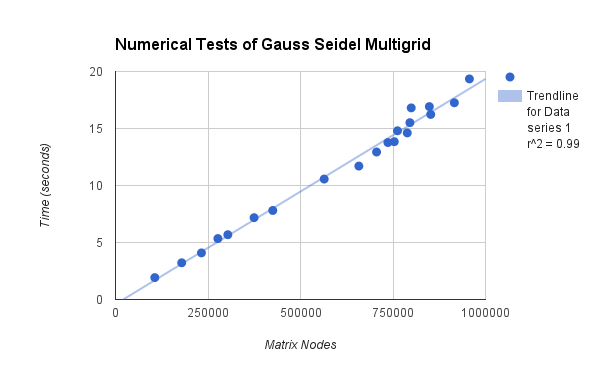}

\section{Conclusion}

In this paper, we have presented an improvement on the existing Cascadic Multigrid Algorithm by introducing a Gauss Seidel smoother as opposed to a power iteration smoother on each level. The algorithm is effective in calculating the algebraic connectivity and the Fiedler vector and is able to partition graphs quickly. 

Having shown that the Gauss-Seidel Cascadic Multigrid Algorithm runs in linear time, we can now discuss its benefits and pitfalls. If our initial graph Laplacian is not sparse, then Gauss Seidel will fail as a smoother since it is inherently meant to work on sparse matrices. In this case, other multigrid algorithms would be optimal. However, the Gauss Seidel smoother works well for most Laplacians as most Laplacians are sparse. Furthermore, we showed that the algorithm is effective in calculating the Fiedler vector of a variety of different graphs. 

We see future works modifying the smoother more. Future improvements could include changing the Gauss Seidel to a Lanczos smoother. Krylov subspace methods are costly for calculating the eigenvalues and eignvectors of large matrices but produce accurate results. Furthermore, future works could include a convergence analysis of the cascadic multigrid algorithm on a more general level and take into account the Gauss Seidel method in the convergence. Of particular interest is our algorithm's convergence with respect to elliptic eigenvalue problems. 

\section{Acknowledgement}

The research presented here was undertaken by Shivam Gandhi and directed by Dr. Xiaozhe Hu of the Tufts University Mathematics Department.

\end{document}